\documentclass[12pt]{amsart}

\usepackage{geometry}               
\usepackage{graphicx}
\usepackage{amssymb}
\usepackage[all]{xy}
\def\aut{\mathop{\text{aut}}}
\def\autn{\mathop{\text{\bf aut}}}
\def\Toro{\mathbb T}
\def\F{\mathop{\mathbb F}}
\def\Z{\mathbb Z}
\def\R{\mathbb R}

\def\alg{\mathop{\text{\bf alg}}}
\def\grp{\text{\bf grp}}
\def\GL{\mathop{\bf GL}}
\def\Diag{\mathop{\bf Diag}}
\def\endo{\mathop{\mathcal{E}nd}}
\def\o{\otimes}
\def\s#1#2{\ {_{#1}}\!\o_{#2}}

\newtheorem{dr}{Drawback}
\newtheorem{te}{Theorem}
\newtheorem{co}{Corollary}
\newtheorem{de}{Definition}
\newtheorem{re}{Remark}
\title{On the gauge action of a Leavitt path algebra}
\author[M. G. Corrales]{Mar\'{\i}a Guadalupe Corrales Garc\'{\i}a}
\address{M. G. Corrales: Centro Regional Universitario de Cocl\'e: \lq\lq Dr. Bernardo Lombardo\rq\rq, Universidad de Panam\'a. Apartado Postal 0229.
Penonom\'e, Provincia de Cocl\'e. Panam\'a.}
\email{mcorrales@ancon.up.ac.pa}
\author[D. Mart\'{\i}n]{Dolores Mart\'{\i}n Barquero}
\address{D. Mart\'{\i}n: Departamento de Matem\'atica Aplicada, Escuela T\'ecnica Superior de Ingenieros Industriales, Universidad de M\'alaga, 29071 M\'alaga, Spain.}
\email{dmartin@uma.es}

\author[C. Mart\'{\i}n]{C\'andido Mart\'{\i}n Gonz\'alez}
\address{C. Mart\'{\i}n: Universidad de M\'alaga, Departamento de \'Algebra Geometr\'{\i}a y Topolog\'{\i}a, Fa\-cultad de Ciencias, Campus de Teatinos s/n, 29071 M\'alaga, Spain.}
\email{candido@apncs.cie.uma.es}

\thanks{The authors have been supported by the Spanish
MEC and Fondos FEDER through project MTM2010-15223, jointly by the Junta
de Andaluc\'{\i}a and Fondos FEDER through projects FQM-336 and FQM-7156.}


\begin{document}
\maketitle
\begin{abstract}
We introduce a revised notion of gauge action in relation with Leavitt path algebras. This notion is based on group
schemes and captures the full information of the grading on the algebra as it is the case of the gauge action of
the graph $C^*$-algebra of the graph.
\end{abstract}



\section{Notations and preliminaries}
For a graph $E$ denote by $C^*(E)$ the graph $C^*$-algebra (see for instance \cite{Rae}) and given any commutative unitary ring $K$, denote by $L_K(E)$ the
Leavitt path algebra associated to $E$ (see \cite[Definition 2.5]{Tom_com} or \cite[p.90]{Workshop} for the case of a ground field of scalars).  
For a graph $E$ we will denote by $E^0$ the set of vertices and $E^1$ the set of edges of $E$. The notation $\hbox{path}(E)$ will be reserved to the set of all path in the graph. As usual given edges $f_1,\ldots, f_n\in E^1$ and the path $\lambda=f_1\ldots f_n$, we will denote $s(\lambda)=s(f_1)$ the source of $f_1$
and $r(\lambda)=r(f_n)$ the range of $f_n$. Also recall that a vertex $v\in E^0$ is said to be regular when $s^{-1}(v)$ is a nonempty finite set. The set of all regular vertices of $E$ will be denoted $\hbox{Reg}(E)$.
In \cite[Corollary 1.5.11]{book},  it  is constructed a basis for $L_K(E)$ when $K$ is a field. Following \cite{book}, the basis can be described taking for each $v\in\hbox{Reg}(E)$, an enumeration of $s^{-1}(v)$ in the form $\{e_1^v,\ldots e_{n_v}^v\}$.
Then, the  basis is given by
$$\{\lambda\mu^*\colon \lambda,\mu\in\hbox{path}(E), r(\lambda)=r(\mu)\}\setminus\{\lambda e_{n_v}^v(e_{n_v}^v)^*\nu^*\colon r(\lambda)=r(\nu)=v\in\hbox{Reg}(E)\}.$$
Observe that the vertices belong to this basis.
As the reader can check, this set is also a basis of $L_K(E)$ in the general case, that is, for any commutative unitary ring $K$. In this way the Leavitt path algebra $L_K(E)$ is always a free $K$-module. Now, some easy observations are convenient for further quotation: 
\begin{re}\label{ffq} \rm
Let $K$ be a commutative unitary ring and $U$, $V$ modules over $K$.
\begin{enumerate}
\item If $U$ is a free $K$-module and $u\in U$ is an element of some basis of $U$ then for any $k\in K$ we have $ku=0$ implies $k=0$.
\item If $U$ and $V$ are free $K$-modules then $U\otimes_K V$ is also a free $K$-module. In particular if $u$ and $v$ are basic elements of $U$ and $V$ respectively, and $r\in K$, the equality $r u\otimes v=0$ in $U\otimes_K V$ implies $r=0$.
\item If $E$ and $F$ are graphs and $(u,v)\in E^0\times F^0$, $r\in K$, then
if $r u\otimes v=0$ in $L_K(E)\otimes L_K(F)$ we have $r=0$.
\end{enumerate}
\end{re}
As usual for any ring $K$ we will denote by $K^\times$ the group of
invertible elements of $K$. Also denote by $\Toro :=S^1$ the unit circle in $\R^2$.
In this work we shall have the occasion to deal with $\Z$-graded algebras. 

The notion of $\Z$-graded algebra $A$ (in a purely algebraic context) is clear: 
the algebra $A$ splits as a direct sum $A=\oplus_{n\in\Z}A_n$ of submodules $A_n$ verifying $A_nA_m\subset A_{n+m}$ for any $n,m\in\Z$.
However, a  $\Z$-grading on a $C^*$-algebra $A$ must be understood 
as defined in \cite[Definition 3.1]{Exel}: $A$ is the closure of a direct sum
$\oplus_{n\in\Z}A_n$ of closed (linear) subspaces $A_n$ of $A$ such that
$A_n^*=A_{-n}$ and $A_nA_m\subset A_{n+m}$ for any $n,m\in\Z$.
We will denote this fact by writing $A=\overline{\oplus_{n\in\Z}A_n}$. 

 When we speak of the canonical $\Z$-grading on $C^*(E)$, we will mean the $\Z$-grading (in $C^*$ sense) such that the component of degree $n$ is formed by those elements $x$ satisfying $\rho(z)(x)=z^n x$ for any $z$.
Roughly speaking, this means that the vertex \lq\lq are of degree $0$\rq\rq\ and the element $f_1\cdots f_ng_1^*\cdots g_m^*$ \lq\lq is homogeneous
of degree $n-m$\rq\rq\ (for any collection of edges $f_i$ and $g_j$).
On the other hand the canonical $\Z$-grading on
$L_K(E)$ is the one for which the vertex are of degree $0$ and the element $f_1\cdots f_ng_1^*\cdots g_m^*$ is homogeneous
of degree $n-m$ (for any collection of edges $f_i$ and $g_j$). Now let  
$A=L_K(E)$ and consider the canonical $\Z$-grading on $A$. Then for any $n$ we consider the canonical
epimorphism $p\colon \Z\to\Z_n$ and the grading on $A$ 
whose component of degree $i$ is the sum $\oplus_{p(n)=i} A_n$. 
This is a coarsening of the canonical $\Z$-grading and since it is a $\Z_n$-grading, we call it the canonical $\Z_n$-grading on $A$. 

\section{Drawbacks of the conventional definition}

The gauge action of the $C^*$-algebra $A:=C^*(E)$
of a graph $E$ is defined as the group homomorphism $\rho\colon \Toro\to\aut(A)$ such that $\rho(z)(p_u)=p_u$ for each
vertex $u$ of the graph and $\rho(z)(s_f)=z s_f$, $\rho(z)(s_f^*)=z^{-1} s_f^*$ for any arrow $f$ and any $z\in \Toro$ (see \cite[Proposition 2.1]{Rae}). With this
definition of the gauge action we can recover the homogeneous components of the canonical $\Z$-grading on $A$ easily, since for any integer
$n$ we have that $A_n$ is just the set
of all $a\in A$ such that for any $z\in \Toro$ we have $\rho(z)(a)=z^n a$. Thus, if we are given the gauge action on $A$, we reconstruct immediately the canonical $\Z$-grading.  Since the gauge action of $A$ codifies all the information of the graded algebra $A$, all the notions
related to this graded structure can be defined in terms of the action. The gauge action is omnipresent in the theory of graph $C^*$-algebras for the same reason that the canonical grading on Leavitt path algebras appears in many of the contributions on the subject. Most of the research works on graph $C^*$-algebras involve its gauge action.  By contrast, most  works on Leavitt path algebras miss the gauge action in the terms in which it has been defined in the literature.

Let us think about the \lq\lq official\rq\rq\
definition of the gauge
action of a Leavitt path algebra $B:=L_K(E)$ over the commutative (and unitary) ring $K$ (see \cite{Ahn}). This is nothing but the group homomorphism $\tau\colon K^\times
\to\aut(B)$ such that $\tau(z)(u)=u$, $\tau(z)(f)=zf$ and $\tau(z)(f^*)=z^{-1}f^*$ for any vertex $u$, any
edge $f$ and any $z\in K^\times$. Though in the original definition, $K$ is a field, we have allowed $K$ to be a unital commutative ring so as to cover the general notion of a Leavitt path algebra. Thus $K^\times$ in the above definition must be understood as  the group of invertible elements of the ring $K$. 
\begin{re}\rm
Let $\tau\colon K^\times\to\aut(A)$ be any representation of the group $K^\times$ by automorphisms on the $K$-algebra $A$ (in particular this applies to the gauge action of a Leavitt path algebra). The action of $\tau$ on an element $t\in K^\times$ will be denoted $\tau(t)$ or $\tau_t$ depending on the typographical convenience.  
\end{re}

Let us  analyze now some peculiarities of the definition above.

\begin{dr}\label{drbone}
The gauge action $\tau$ does not capture the whole information of the graded algebra $L_K(E)$.
\end{dr}
Proof. In fact in some extreme cases $\tau$  contains no information at all simply because $\tau$ is trivial. For instance
take $K=\F_2$ to be the field of two elements. Then $K^\times$ is the trivial group $K^\times=\{1\}$ and $\tau$ is
the trivial group homomorphism $1\mapsto 1$. So in this case $\tau$ gives no information at all of the grading on $B=L_K(E)$. 
In other cases in which $\tau$ is not trivial, we can not recover the original grading on $B$. Take for instance $K$ to be any ring whose group of invertibles
$K^\times$ is isomorphic to $\F_2$ (for instance $K$ could be the ring $\Z_4$ of integers module $4$, or even the field of three elements $K=\F_3$).
Since $K^\times=\{1,\omega\}$ with $\omega^2=1$, then $\tau$ is completely determined by $\tau(\omega)$. If we denote by $B_n$
the homogeneous component of degree $n$ of the canonical grading on $B$, then  $B_n$ does not agree with
the submodule
$$\{x\in B\colon \tau(z)(x)=z^n x, \text{for all } z\in K^\times\}.$$ 
In fact the above submodule agrees with $\{x\in B\colon \tau(\omega)(x)=\omega^n x\}$ and agrees with
the direct sum of homogeneous components of even degree (if $n$ is even) and the sum of components of odd degree if $n$ is odd.  So what we get from the gauge action is the canonical $\Z_2$-grading obtained as a coarsening of the
canonical $\Z$-grading of $B$. By contrast with what happens to the gauge action on $C^*(E)$ we have
\begin{equation}\label{incl}
B_n\subset \{x\in B\colon \tau(z)(x)=z^n x, \text{for all } z\in K^\times \}
\end{equation}
but equality rarely holds when $K$ is a ring (or even for some fields).
It can be proved that for an infinite field the inclusion in (\ref{incl}) is in fact an equality. Thus the \lq\lq official\rq\rq\  definition only works adequately on infinite fields. Our proposal is to give a definition working correctly on any ring of scalars.  

\begin{dr}\label{drtwo} For the gauge action $\rho$ of $A=C^*(E)$ the notion of graded ideal is equivalent to that of $\rho$-invariant ideal. This is not the case for the gauge action $\tau$ of $B=L_K(E)$.
\end{dr}
Proof. Recall that a closed ideal $I$ of $A$ is $\rho$-invariant if and only if 
$\rho(z)(I)\subset I$ for all $z\in \Toro$. Of course when
$I$ is a graded ideal then it is $\rho$-invariant. On the other hand if $\rho(z)(I)\subset I$ for any $z\in \Toro$ it is easily seen that $I$ is graded: take $x\in I$ and consider its expansion as a series of homogeneous elements $x=\sum_n x_n$.
We must prove that each $x_n$ is an element in $I$. We know that
$x_n=\int_{S^1} z^{-(n+1)}\rho(z)(x)dz$   (normalizing the integral on $S^1$ so that $\int_{S^{1}}\frac{dz}{z}=1$). Thus taking into account that $\rho(z)(I)\subset I$ and that $I$ is closed, we get $x_n\in I$.

On the other hand if $\tau$ denotes the gauge action of $B$ and we take a $\tau$-invariant ideal $I$ of $B$, we can not prove $I$ to be graded. 
In an attempt to do it,  taking $a\in I$ we know that
$a=\sum_n a_n$ where $a_n\in B_n$ (homogeneous component of degree $n$) then  $\tau(z)(a)\in I$ hence
$\sum_n z^n a_n\in I$ for any $z\in K^\times$. Thus 
\begin{equation}\label{zuru}
(a_{-k},\ldots, a_0,a_1,\ldots, a_k)
\begin{pmatrix} 1 & z_1^{-k} & \cdots & z_{2k}^{-k}\cr
1 & z_1^{-k+1} & \cdots & z_{2k}^{-k+1}\cr
\vdots & \vdots & \ddots             & \vdots\cr 
1 & z_1^{k} & \cdots & z_{2k}^{k}
\end{pmatrix}\in\overbrace{ I\times\cdots\times I}^{2k+1}=:I^{2k+1}\end{equation}
for any collection of $z_i\in K^\times$. Denoting by $M$ the matrix in the left hand member of formula (\ref{zuru}) if this matrix were invertible
then could conclude that $$(a_{-k},\ldots, a_0,a_1,\ldots, a_k)\in I^{2k+1}$$ hence we would have $a_n\in I$ for every $n$. But
the determinant of $M$ is a Vandermonde determinant and it would have to be an invertible element in $K$ in order to have $M$ invertible. Not even if $K$ is a field can we say for sure that this determinant is nonzero. The $z_i$s would have to be all different but this collection of scalars may be larger than the cardinal of $K$. 
So for the Leavitt path algebra $B$ the possibility of choosing these scalars all different depends on the ground ring of scalars, and it is not always guaranteed. 
Thus, in general, we do not have that gauge invariant ideals of $B$ are graded ideals.

\begin{dr}\label{drthree} Consider two graph $C^*$-algebras $A_i$ ($i=1,2$) with associated gauge actions $\rho_i$. Define a homomorphism $f\colon A_1\to A_2$ to be a gauge-homomorphism when for any $z\in \Toro$ the following square is commutative
\[
\xymatrix{
A_1 \ar[r]^f\ar[d]_{\rho_1(z)} & A_2\ar[d]^{\rho_2(z)}\\
A_1 \ar[r]_f & A_2
}
\]
In a similar fashion can we define the notion of a gauge homomorphism of Leavitt path algebras. However while
the notion of gauge-homomorphism is equivalent to that of graded homomorphism in the setting of graph $C^*$-algebras, it is not the case that for Leavitt path algebras, both notions agree.
\end{dr}
Proof. In the ambient of Leavitt path algebras the gauge action may be even trivial (if the ground field has characteristic $2$). Thus gauge-homomorphisms are simply homomorphism in this case. Since not every homomorphism is graded
we see that both notions do not agree. For graph $C^*$-algebras it is easily seen that both notions agree. Of course this drawback and the previous do not exist for Leavitt path algebras over infinite fields
but we would like to give a notion of gauge action which overcomes these difficulties and does not
depend so much of the ground ring of scalars.
\medskip

Some more drawbacks will be explained in the sequel. For the moment we can realize that the canonical $\Z$-grading of graph $C^*$-algebras and also in Leavitt path algebras are present in many argumentations on these algebras. 
Frequently we can argue on homogeneous elements and then generalize to arbitrary ones. That is why is so frequent
to find the gauge action in the graph $C^*$-algebra literature and not in the Leavitt path algebra one (where one must replace gauge action arguments with others involving the canonical $\Z$-grading).  

\begin{dr} 
The Gauge-Invariant Uniquenes Theorem is stated in \cite{Rae} in the following terms
\begin{te} {\rm (\cite[Theorem 2.2, p.16]{Rae})}
Let $E$ be a row-finite graph and suppose that $\{T,Q\}$ is a Cuntz-Krieger $E$-family in a $C^*$-algebra $B$
with each $Q_v\ne 0$. If there is a continuous action $\beta\colon \Toro\to\aut(B)$ such that $\beta(z)(T_e)=zT_e$
for every $e\in E^1$ and $\beta(z)(Q_v)=Q_v$ for every $v\in E^0$, then $\pi_{T,Q}$ is an isomorphism onto
$C^*(T,Q)$.
\end{te}

As far as we know the best version of the previous theorem for Leavitt path algebras is given in \cite[Theorem 1.8, p. 6]{Ahn} and it claims:

\begin{te}\label{AGIUTprimi}
(The Algebraic Gauge-Invariant Uniqueness Theorem.)  Let $E$ be a row-finite graph,  $K$ an
infinite field, and $A$ a $K$-algebra. Denote by $\tau^E$ the gauge action of $L_K(E)$. Suppose
$$\phi\colon L_K(E)\rightarrow A$$
 is a $K$-algebra homomorphism such that $\phi (v)\ne 0$ for every $v\in E^0$. If there exists a group
action $\sigma \colon K^\times \rightarrow \mbox{Aut}_K(A)$ such that $\phi \circ \tau^E _t=\sigma _t \circ
\phi$ for every $t\in K^\times$,  then $\phi $ is injective.
\end{te}
The hypothesis on the infiniteness of the ground field can not be removed as the following example shows: take $K=\Z_2$ and $A=L_K(E)/I$ where $I$ is an ideal which does not contain any vertex. For instance $E$ could be the one-petal rose, (where $E^0$ and $E^1$ have cardinal $1$). Then $L_K(E)\cong K[T,T^{-1}]$ and the ideal $I$ generated by the non-invertibe element $1+T$ does not contain any vertex. So with these ingredients the canonical epimorphism $p\colon L_K(E)\to A$ satisfies $p(v)\ne 0$ for (the unique) $v\in E^0$ and it is not injective. Furthermore
the gauge action of $L_K(E)$ is trivial (since $K^\times=\{1\}$) and we can
consider the trivial group action $\sigma\colon K^\times\to\aut(A)$ and the requirement $p\circ\tau_t^E=\sigma_t\circ p$ is trivially satisfied.

We can give another  example in which $K^\times$ is not trivial.
So assume $K$ to be a ring such that $K^\times\cong\F_2$.
Then $K^\times=\{1,\omega\}$ with $\omega^2=1$. Consider as before
$E$ the one-petal rose and $A=L_K(E)/I$ where $I$ is the ideal generated
by the noninvertible element $1+T^2$ (so $I$ does not contain the unique vertex
of $E^0$). The gauge action $\tau^E$ is completely defined by $\tau^E_\omega$
and we have $$\tau^E_\omega(T^n)=
\begin{cases} T^n \hbox{ if } n \hbox{ is even}\cr
\omega T^n \hbox{ otherwise.}\end{cases}$$
Next we define $\sigma\colon K^\times\to\aut(A)$ by 
$$\sigma_\omega(T^n+I):=
\begin{cases} T^n+I \hbox{ if } n \hbox{ is even}\cr
\omega T^n +I\hbox{ otherwise.}\end{cases}$$
As before, the canonical epimorphism $p\colon L_K(E)\to A$ satisfies
$p(v)\ne 0$ for  the unique $v\in E^0$ and $p\circ\tau^E_t=\sigma_t\circ p$ (for $t\in K^\times$) but $p$ is not injective.
\smallskip

As we shall see, the hypothesis on the infiniteness of $K$ in Theorem \ref{AGIUTprimi} is not necessary if we use the schematic version of
the gauge action. 
\end{dr}

\begin{dr}
The natural translation of the crossed product of $C^*$-algebras to a purely algebraic setting must be made carefully.
\end{dr}
We recall the definition of the crossed product of $C^*$-algebras. Assume that $A$ and $B$ are $C^*$-algebras and $G$ a compact abelian group with actions $\mu\colon G\to\aut(A)$ and $\nu\colon G\to\aut(B)$. Consider next the action
$\lambda\colon G\to \aut(A\otimes B)$ defined by $\lambda(g)(a\otimes b)=\mu(g)(a)\otimes \nu(g^{-1})(b)$. Define
now the {\em crossed product} $A\otimes_G B$ as the fixed point algebra under the action $\lambda$. 

The gauge action of a graph $C^*$-algebra has been successfully applied to certain interesting constructions in
\cite{Bates} and \cite{Kunijan}. Take two row-finite graphs $E$ and $F$ and define its product $E\times F:=(E^0\times F^0,E^1\times F^1,r,s)$ where $s(f,g)=(s(f),s(g))$ and $r(f,g)=(r(f),r(g))$ for any $(f,g)\in E^1\times F^1$.  Though this is not the usual
definition of the product of two graphs, this notion is interesting for us since it allows to describe the crossed
product of graph $C^*$-algebras.  Indeed, it is proved  in \cite[Proposition 4.1, p. 62]{Bates} that if $E$ and $F$ are row-finite graphs with no sinks, then there is an  isomorphism 
\begin{equation}\label{simi}
C^*(E\times F)\cong C^*(E)\otimes_{\Toro}C^*(F)
\end{equation}

\noindent where the crossed product on the right-hand side is the induced by the gauge actions $\Toro\to \aut(C^*(E))$ and $\Toro\to \aut(C^*(F))$.
We can try to mimic this definition of crossed product in a purely algebraic context. So assume that $G$ is an abelian group action by automorphisms in the $K$-algebras $A$ and $B$. Let $\mu\colon G\to\aut(A)$ and $\tau\colon G\to\aut(B)$ be two representations of $G$ in $A$ and $B$ respectively and we define $\lambda\colon G\to\aut(A\otimes B)$ as above. Then the fixed point algebra under the action $\lambda$ will be denoted by $A\otimes_G B$. Thus, $A\otimes_G B$ is a subalgebra of the usual tensor product algebra $A\otimes_K B$. Observe that when the group $G$ is trivial we get $A\otimes_G B=A\otimes_K B$. 

If we use this naive interpretation of the crossed product of algebras, the similar property to (\ref{simi}) for Leavitt path algebras does not hold.
For instance if the consider $L_K(E)$ for a ring such that $K^\times$ is trivial.
Then $L_K(E)\otimes_{K^\times}L_K(E)=L_K(E)\otimes_K L_K(F)$ and we would have $L_K(E\times F)\cong L_K(E)\otimes_K L_K(F)$. But there are two ways to see that this is not true:
\begin{enumerate}
\item Apply the results in \cite{Ara} in which the impossibility of this isomorphism
is studied. More concretely assume $E$ is the one-petal rose ($\vert E^0\vert=\vert E^1\vert=1$) and $F$ is the two-petals rose ($\vert F^0\vert=1$, $\vert F^1\vert=2$). Then $E\times F\cong F$ (a graph isomorphism) and applying
\cite[Theorem 5.1, p.2635]{Ara} for $n=2$, $E_1=E_2=E$, $m=1$ $F_1=F$ we conclude that $L_K(E)\otimes L_K(F)$ is not isomorphic to $L_K(F)$.
  
\item A simple example also proves the impossibility of the isomorphism
$L_K(E\times F)\cong L_K(E)\otimes_K L_K(F)$. You can get the graph $E$ on the left hand of:
 
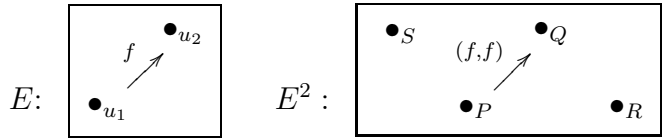
\[
 \begin{tabular}{ccccc}
 $E$: &  \fbox{\xygraph{!{<0cm,0cm>;<1cm,0cm>:<0cm,1cm>::}
 !{(0,0) }*+{\bullet_{u_1}}="a" 
 !{(1,1) }*+{\bullet_{u_2}}="b"
 "a":^f"b" }}
  & & $E^2:$ &
\fbox{ \xygraph{!{<0cm,0cm>;<1cm,0cm>:<0cm,1cm>::}
 !{(0,0) }*+{\bullet_{P}}="a" 
 !{(1,1) }*+{\bullet_{Q}}="b"
 !{(2,0) }*+{\bullet_{R}}="c"
 !{(-1,1) }*+{\bullet_{S}}="d"
 "a":^{(f,f)}"b" "c" "d" }}
 \end{tabular}
 \]
then $E^2:=E\times E$ is the graph on the right hand of the above figure. Thus $L_K(E)\cong \mathcal M_2(K)$ and
$L_K(E^2)=K\oplus K \oplus {\mathcal M}_2(K)$ which has dimension $6$. However $L_K(E)\otimes L_K(E)$ has dimension $16$. Hence the cross product
of the Leavitt path algebras does not agree with the Leavitt path algebra of $E^2$. 
\end{enumerate}

Once we have realized some handicaps of the gauge action of a Leavitt path algebras, we propose a different approach.

\section{Redefining the gauge action}

For any associative and commutative ring $K$ (with unit $1\in K$) denote by  
$\alg_K$ the category of unital associative commutative $K$-algebras and unital homomorphisms. Denote
by $\grp$ the category of groups. Recall that a $K$-group functor is a covariant
functor $\mathcal F\colon \alg_K\to\grp$. If $\mathcal F$ and $\mathcal G$
are $K$-group functors, a homomorphism $\eta\colon{\mathcal F}\to{\mathcal G}$ is nothing but a natural transformation from $\mathcal F$ to $\mathcal G$. 

Recall also that an affine $K$-group scheme is a representable  $K$-group  functor, that is, $\mathcal F=\hom_{\alg_K}(H,\--)$ for some Hopf algebra $H$ (see \cite{Water} or \cite{Jantzen}). An affine group scheme is said to be an algebraic group if the representing Hopf algebra is finitely-generated.

For any $K$-algebra $A$ (not necessarily associative or commutative or unital) we can consider the $K$-group functor $\autn(A)\colon\alg_K\to\grp$
such that $\autn(A)(R):=\aut(A_R)$ (where $A_R:=A\otimes R$) for any object $R$ in the category $\alg_K$. We emphasize that $\aut(A_R)$ denotes the group of $R$-algebra automorphisms of $A_R$. We also recall the definition of the $K$-group functor $\GL_n
\colon\alg_K\to\grp$ such that $\GL_n(R)$ is the group of invertible $n\times n$ matrices with entries in $R$. In particular $\GL_1(R)=R^\times$ the group of
invertible elements in $R$. This $K$-group $\GL_1$ is representable (its representing Hopf algebra being the Laurent polynomial algebra
$K[x,x^{-1}]$), hence it is an affine group scheme (and even an algebraic group).  

A diagonalizable affine group scheme (diagonalizable group in the sequel)
is an affine group scheme whose representing Hopf algebra is the group algebra of an abelian group.  Thus if $\Lambda$ is an abelian group and we consider
the group algebra $K\Lambda$ (with its natural structure of Hopf algebra), then the
$K$-group functor $\hom_{\alg_K}(K\Lambda,\--)$ is said to be diagonalizable and 
its usual notation is $\Diag(\Lambda):=\hom_{\alg_K}(K\Lambda,\--)$. When $\Lambda$ is finitely-generated $\Diag(\Lambda)$ is an algebraic group (the
group structure in $\hom_{\alg_K}(K\Lambda,R)$ is point-wise multiplication, that is, if $\alpha,\beta\in\hom(K\Lambda,R)$ then $(\alpha\beta)(x):=\alpha(x)\beta(x)$ for any $x\in K\Lambda$).

\subsection{Representation of diagonalizable groups}

In this section we note a series of results which are well-known but not so
easy to quote (at least in its present form).  There are two (equivalent) approaches to the study of gradings. Both are based upon affine group schemes. On the one hand we have the co-modules approach which
skips the most puzzling notion of affine schemes, and on the second hand
we have the representations of diagonalizable groups. The gauge action
of $C^*$-algebras as well as the definition of gauge action for Leavitt path
algebras are closer to the  viewpoint of representations of diagonalizable
group schemes. So we adopt this philosophy.  

 Most of the materials in this
subsection can be seen in \cite{Gabriel} and in \cite{Jantzen}. Also in  \cite{Susan}  with the slightly different terminology of co-modules.  That is why we include this subsection
in which we simply translate the main results to the language of representations. Of course the reader familiarized with representations
of affine group schemes could skip this subsection and proceed with
the next.

 Consider a $K$-module $M$ and define the $K$-group functor $\GL(M)\colon\alg_K\to\grp$ such that for any $R$ we have $\GL(M)(R):=\text{GL}(M_R)$ where
$M_R:=M\otimes R$ and $\text{GL}(M_R)$ is the group of invertible automorphisms of the $R$-module $M_R$. 
If $G$ is a $K$-group functor and $\rho\colon G\to\GL(M)$
is a $K$-group homomorphism then it is said that $\rho$ is a representation
of $G$. This admits in a standard way a formulation in terms of modules as the reader can guess. For any $K$-algebra $R$ in $\alg_K$ we have a
group homomorphism $\rho_R\colon G(R)\to \text{GL}(M_R)$ and if $\alpha\colon R\to S$ is a homomorphism of $K$-algebras there is a commutative
diagram
\[\xymatrix{ G(R)\ar[r]^{\rho_R}\ar[d]_{G(\alpha)} &\text{GL}(M_R)\ar[d]^{\alpha^*}\\
G(S)\ar[r]_{\rho_S} & \text{GL}(M_S)}
\]
where for any $f\in\text{GL}(M_R)$ the map $\alpha^*(f)\colon M_S\to M_S$
is given by $\alpha^*(f)(m\otimes 1)=(1\otimes\alpha)f(m\otimes 1)$. In other
words, $\alpha^*$ is the morphism-function of the functor $\GL(M)$.

If $G$ turns out to be a diagonalizable group $G=\Diag(\Lambda)$ and we
consider a representation $\rho\colon G\to\GL(M)$ then we have group
homomorphisms 
$\rho_R\colon \hom(K\Lambda,R)\to\text{GL}(M_R)$
for any $K$-algebra $R$. In particular, we have 
$\rho_{K\Lambda}\colon \endo(K\Lambda)\to\text{GL}(M_{K\Lambda})$
and the above commutative diagram specializes to 
 \[\xymatrix{f\ar@{|->}[d] & \endo(K\Lambda)\ar[r]^{\rho_{K\Lambda}}\ar[d] &
 \text{GL}(M_{K \Lambda})\ar[d]^{\alpha^*}\\
\alpha\circ f& \hom(K\Lambda, R)\ar[r]_{\rho_R} & \text{GL}(M_R)}
\]
where $\alpha\colon K\Lambda\to R$ is a $K$-algebra homomorphism. Then
the commutativity of the diagram yields (taking $f=1_{K\Lambda}$) the formula
 $\rho_R(\alpha)=\alpha^*(\rho_{K\Lambda}(1_{K\Lambda}))$ and so 
 \begin{equation}\label{mg1}
 \rho_R(\alpha)(m\otimes 1)=(1\otimes\alpha)[
 \rho_{K\Lambda}(1_{K\Lambda})(m\otimes 1)]
 \end{equation}
 that is, $\rho_R$ is completely determined by $\rho_{K\Lambda}$ and so
 the whole action $\rho$ is completely determined once we know $\rho_{K\Lambda}$.  If we take $m\in M$ we may write
 $\rho_{K\Lambda}(1_{K\Lambda})(m\otimes 1)=
 \sum_{\lambda\in\Lambda} p_\lambda(m)\otimes \lambda$ where
 $p_\lambda\colon M\to M$ is a $K$-modules homomorphism.
\begin{te}
The set $\{p_\lambda\}_{\lambda\in\Lambda}$ is a system of orthogonal
idempotents of $\endo_K(M)$  and for any $m\in M$ we have
$m=\sum_\lambda p_\lambda(m)$.
\end{te}
Proof. Consider the identity $c$ of the group $\endo(K\Lambda)$. This
is the map such that $c(\lambda)=1$ for any $\lambda\in\Lambda$.  So $\rho_{K\Lambda}(c)$
is the identity on $M_{K\Lambda}$ therefore $m\otimes 1=
\rho_{K\Lambda}(c)(m\otimes 1)$ and applying formula 
(\ref{mg1}) we get 
$$m\otimes 1=
\rho_{K\Lambda}(c)(m\otimes 1)=(1\otimes c)[\rho_{K\Lambda}(1_{K\Lambda})
(m\otimes 1)]=$$ $$(1\otimes c)(\sum_\lambda p_\lambda(m)\otimes \lambda)=
\sum_\lambda p_\lambda(m)\otimes 1$$ whence the second assertion of the Theorem ($m\otimes 1=0$ implies $m=0$ applying $1\o\epsilon$ where
$\epsilon$ is the counit of the Hopf algebra $K\Lambda$). To see the first one we use: the equality 
$\rho_{R}(\alpha\beta)=\rho_{R}(\alpha)\rho_{R}(\beta)$ which holds for any
$K$-algebra $R$ and any
 $\alpha,\beta\in\hom_{\alg_K}(K\Lambda,R)$. Thus we take $R=K\Lambda\o K\Lambda$ and $\alpha,\beta\colon K\Lambda\to R$ such that $\alpha(\lambda)=\lambda\o 1$, $\beta(\lambda)=1\o\lambda$ for any $\lambda\in\Lambda$.
 
  Taking into account this as well as equation (\ref{mg1}) we get:
$$\sum_\lambda p_\lambda(m)\otimes\lambda\o\lambda=
\sum_{\lambda,\mu}p_\mu p_\lambda(m)\otimes \mu\o\lambda$$
and since $\{\mu\o\lambda\colon \mu,\lambda\in \Lambda\}$ is a basis of $R$
we conclude that for any $\lambda,\mu\in\Lambda$ one has $p_\lambda^2=p_\lambda$ and $p_\mu p_\lambda=0$ if
$\lambda\ne \mu$.

\begin{co}
If the $K$-module $M$ admits a decomposition 
$M=\oplus_{\lambda\in\Lambda} M_\lambda$ where $\Lambda$ is an abelian group,
then there is a representation $\rho\colon\Diag(\Lambda)\to
\GL(M)$ such that for any $K$-algebra $R$ the map $\rho_R$ acts in
the form $\rho_R(\alpha)(m_\lambda\otimes 1)=m_\lambda\otimes\alpha(\lambda)$ for any $\alpha\in\hom(K\Lambda,R)$ and any
$m_\lambda\in M_\lambda$. Reciprocally given a representation 
$\rho\colon\Diag(\Lambda)\to
\GL(M)$ of the diagonalizable affine group scheme $\Diag(\Lambda)$,
there is a decomposition $M=\oplus_\lambda M_\lambda$ where
$M_\lambda=p_\lambda(M)$.
\end{co}
Proof. The only thing to take into account for the reciprocal is that the
set of orthogonal idempotents $\{p_\lambda\}$ induce the decomposition
 $M=\oplus_\lambda M_\lambda$.
 
 Now the last result can be adapted to handle group gradings on algebras.
  \begin{co}\label{kaka}
If $A$ is a $K$-graded algebra $A=\oplus_{\lambda\in\Lambda}A_\lambda$  where $\Lambda$ is an abelian group,
there is a representation $\rho\colon\Diag(\Lambda)\to
\autn(A)$ such that for any $K$-algebra $R$ the map $\rho_R$ acts in
the form $\rho_R(\alpha)(x_\lambda\otimes 1)=x_\lambda\otimes\alpha(\lambda)$ for any $\alpha\in\hom(K\Lambda,R)$ and any
$x_\lambda\in A_\lambda$. Reciprocally given a representation 
$\rho\colon\Diag(\Lambda)\to
\autn(A)$ of the diagonalizable affine group scheme $\Diag(\Lambda)$,
there is a grading $A=\oplus_{\lambda\in\Lambda}A_\lambda$ where
$A_\lambda$ is the set of all $x\in A$ such that $\rho_R(\alpha)(x\otimes 1)=
x\otimes\alpha(\lambda)$ (for any $K$-algebra $R$ and $\alpha\in\hom(K\Lambda,R)$).
\end{co}

\subsection{A new concept of gauge action}

For a Leavitt path algebra $A=L_K(E)$ over a ring $K$, the gauge action is defined as the 
map $\rho\colon K^\times\to\aut(A)$ such that for any $z\in K^\times$
the automorphism $\rho(z)$ fixes the vertices and $\rho(z)(f)=zf$, $\rho(z)(f^*)=z^{-1}f^*$ for any edge $f$. The definition we propose implies affine groups schemes. So consider the affine group scheme $\GL_1=\Diag(\Z)$ whose
representing Hopf algebra is the group algebra $K\Z$ of $\Z$ which we identify
with the Laurent polynomial algebra $K[x,x^{-1}]$. Thus $\GL_1:=\hom(K[x,x^{-1}],\--)$ and we have $\GL_1(R)=\hom(K[x,x^{-1}],R)\cong R^\times$ since a $K$-algebras homomorphism $K[x,x^{-1}]\to R$ is completely
determined by the image of $x$ in $R$ (which is an invertible element in $R$, and so it belongs to $R^\times$). Consequently we will identify $\GL_1(R)$ with the
group of invertible elements $R^\times$ of $R$. Obviously $\GL_1(K)=K^\times$. 
\begin{de}
For a Leavitt path algebra $A=L_K(E)$ over a ring $K$ define the
gauge action as a representation of the diagonalizable group scheme
$\GL_1$ given by $\rho\colon\GL_1\to\autn(A)$ where for any $K$-algebra
$R$ and any $z\in R^\times$ we have $\rho_R(z)(u\otimes 1)=u\otimes 1$ for
each vertex $u$ and 
$\rho_R(z)(f\otimes 1)=f\otimes z$,  $\rho_R(z)(f^*\otimes 1)=f^*\otimes z^{-1}$
for any $f\in E^1$.
\end{de}

Whith this schematic approach we see that $\rho_K$ agrees with the
\lq\lq official\rq\rq\ definition of gauge action. On the other hand we can define such an action for any $K$-algebra $A$ endowed with a $\Z$-grading
$A=\oplus_{n\in\Z}A_n$: just define for any $K$-algebra $R$ the map
$\rho_R(z)(a_n\otimes 1):=a_n\otimes z^n$ for any $n\in\Z$, $a_n\in A_n$
and any $z\in R^\times$.
Reciprocally if we have a representation
$\rho\colon\GL_1\to\autn(A)$ for some $K$-algebra $A$, then applying Corollary \ref{kaka}, $A$ is $\Z$-graded where for each integer $n$ we have 
\begin{equation}\label{adidaone}
A_n=\{a\in A\colon
\rho_R(z)(a\otimes 1)=a\otimes z^n \text{ for all } z\in R^\times \text{ and each object } R \text{ in } \alg\nolimits_K\}.
\end{equation}

The first conclusion we get is 
\begin{te}
The gauge action in schematic sense encloses all the information of the grading. We can recover the homogeneous components from the schematic gauge action. So Drawback \ref{drbone}
no longer holds with this new definition. 
\end{te}

\begin{re}\rm
At this point we can see how the new gauge action fixes the problem posed
in Drawback \ref{drbone}, which affects to $B:=L_K(E)$ when the ring $K$ has a trivial group $K^\times=\{1\}$. Of course in this case the \lq\lq official\rq\rq\  gauge action is trivial. By using the new notion we have for any (associative, commutative and unital ) $K$-algebra $R$, a representation $\rho_R\colon R^\times\to\aut(B\otimes_K R)$. Furthermore $B_n$ agrees with the set of all
$x\in B$ such that $\rho_R(z)(x\otimes 1)=x\otimes z^n$ for any $z\in R^\times=\GL_1(R)$.
Take for instance the group algebra $R:=K\Z$. This can be identified with the
Laurent polynomial $K$-algebra $K[T,T^{-1}]$ for some indeterminate $T$.
After that identification,  $R^\times$ contains the set of elements $\{T^k\colon k\in\Z\}$ (which is infinite independently of the nature of $K$). Then the particular 
representation $\rho_{K\Z}\colon (K\Z)^\times\to\aut(B\otimes K\Z)$ suffices to
describe the homogeneous components $B_n$, ($n\in\Z$).  Indeed, by (\ref{adidaone}) we have
$$B_n=
\{x\in B\colon \rho_{K\Z}(z)(x\otimes 1)=x\otimes z^n, z\in (K\Z)^\times\}.$$
But even the weaker set of conditions $\rho_{K\Z}(T^k)(x\otimes 1)=
x\otimes T^{kn}$ for any $k\in\Z$, also implies $x\in B_n$: decompose
$x=\sum x_q$ with $x_q\in B_q$; then $\sum_q \rho_{K\Z}(T^k)(x_q\otimes 1)=\sum_q x_q\otimes T^{kq}$. So $\sum_q x_q\otimes T^{kq}= \sum_q x_q\otimes T^{kn}$ and given the linear independence of the powers of $T$ we conclude
$x_q=0$ for any $q\ne n$. Thus $x\in B_n$. 
\end{re}

\begin{re}\rm 
The key point of the proposed new definition, it that it allows us to extend scalars to any $K$-algebra and to recover the information on the original algebra. So it is a going-up and going-down process. Of course this new notion gives more information (retains all the information of the canonical grading which is not guaranteed  by the old notion) because it is also more demanding.
The scalar extension procedure that we propose forms part of the essentials of
group scheme theory which is what we are applying here.  
\end{re}

Let us go now to the notion of $\rho$-invariant ideal of $A:=L_K(E)$.
So we assume given $\rho\colon\GL_1\to\autn(A)$ the gauge action
in schematic sense. 
\begin{de}
An ideal $I$ of $A$ is said to be $\rho$-invariant when 
for any $K$-algebra $R$ and any $z\in R^\times$ we have
$\rho_R(z)(I\otimes 1)\subset I\otimes R$.
\end{de}
Clearly, if $I$ is a graded ideal of $A$ then $I$ is $\rho$-invariant: indeed
take a $K$-algebra $R$ and any $z\in R^\times$. Take any $a\in I$, then
$a=\sum a_n$ where each $a_n\in I\cap A_n$. Thus
$\rho_R(z)(a\otimes 1)=\sum_n \rho_R(z)(a_n\otimes 1)=\sum_n a_n\otimes
z^n\in I\otimes R$. Consequently graded ideals of $A$ are $\rho$-invariant.
But the reciprocal is also true:

\begin{te}\label{thrive}
An ideal $I$ of $A$ is graded if and only if it is $\rho$-invariant. Thus
Drawback \ref{drtwo} no longer holds.
\end{te}
Proof.
Let $R:=K[T,T^{-1}]$ be the Laurent polynomial algebra in the indeterminate $T$ over the commutative unitary ring $K$.
Then $\{T^n\colon n\in\mathbb{Z}\}$ is a linearly independent set.
Define the $K$-modules homomorphism $f_n\colon R\to K$ by $f_n(T^m)=\delta_{nm}$ (Kronecker's delta). Consider the $K$-bilinear map
$A\times R\to A$ such that $(a,r)\mapsto f_n(r)a$ and the $K$-modules homomorphism $\Phi_n\colon A\otimes R\to A$
such that $\Phi_n(a\otimes r)=f_n(r)a$. If $I$ is an ideal in $A$ then $\Phi_n(I\otimes R)\subset I$.
If $\rho\colon\mathbf{GL}_1\to\mathbf{aut}(A)$ is the gauge action,
take $a\in I$ and decompose it as  $a=\sum_m a_m$ where $a_m\in A_m$.
Assume that $I$ is $\rho$-invariant. Then $\rho_R(T)(a)=\sum_m a_m\otimes T^m\in I\otimes R$.
So $\Phi_n(\sum_m a_m\otimes T^m)\in I$, but then $I\ni\Phi_n(\sum_m a_m\otimes T^m)=\sum_m f_n(T^m)a_m=a_n$ for any $n$. Thus $I$ is graded.
\medskip

Let us deal with Drawback \# \ref{drthree} now.
Given two Leavitt path $K$-algebras $A_1$ and $A_2$ with their respective gauge actions  in schematic sense $\rho_i\colon\GL_1\to\autn(A_i)$, $i=1,2$.
Then
\begin{de}\label{cuad}
A homomorphism $f\colon A_1\to A_2$ is said to be a gauge-homomorphism if
the following square is commutative.
\[
\xymatrix{
A_1\o R \ar[r]^{f\o 1}\ar[d]_{(\rho_1)_R(z)} & A_2\o R\ar[d]^{(\rho_2)_R(z)}\\
A_1\o R \ar[r]_{f\o 1} & A_2\o R
}
\]
for any $K$-algebra $R$ and any $z\in R^\times$.
\end{de}
It is easy to prove that in case $f\colon A_1\to A_2$ is a graded homomorphism, then it is a gauge-homomorphism: take $a\in A_1$ homogeneous of degree
say $n$. Then $(\rho_2)_R(z) (f\o 1)(a\o 1)=(\rho_2)_R(z) (f(a)\o 1)$ and since
$f(a)$ is an homogeneous element of $A_2$ of degree $n$ then $(\rho_2)_R(z) (f\o 1)(a\o 1)=f(a)\o z^n=(f\o 1)(a\o z^n)=
(f\o 1)(\rho_1)_R(z)(a\o 1)$. Thus $f$ is a gauge-homomorphism.
But we have also the reciprocal.
\begin{te}\label{gradone}
The homomorphism $f\colon A_1\to A_2$ is graded if and only if it is a gauge-hom\-omor\-phism. Thus Drawback \ref{drthree} disappears.
\end{te}
Proof.
Assume that $f$ is a gauge-homomorphism and take $a$ in the homogeneous component of degree $n$ of $A_1$.  Recall such component agrees with the submodule of all the elements $a\in A_1$ such that $(\rho_1)_R(z)(a\otimes 1)=a\otimes z^n$ for any (unitary associative and commutative) $K$-algebra $R$ and any $z\in R^\times$. Consequently we must prove that $(\rho_2)_R(z)(f(a)\otimes 1)=a\otimes z^n$ for any $R$ and $z$ as before. But this is a direct corollary of the commutativity of the square in Definition \ref{cuad}.
\medskip

\begin{re}\label{notilla}\rm
Theorem \ref{gradone} can be generalized in the following sense: let $A_i$ ($i=1,2$) be $K$-algebras endowed with representations $\rho_i\colon\GL_1\to\autn(A_i)$, ($i=1,2$). Consider now the $\Z$-grading induced by $\rho_i$ in $A_i$: the homogeneous component of degree $n$ is just the $K$-submodule of those elements $a\in A_i$ such that $(\rho_i)_R(z)(a\otimes 1)=a\otimes z^n$ for any associative, commutative and unitary $K$-algebra $R$, and any $z\in R^\times$.   Take then a $K$-algebra homomorphism $f\colon A_1\to A_2$. With the same proof as above, we have that  $f$ is graded if and only if it is a gauge-homomorphism in the sense that the
squares in Definition \ref{cuad} are commutative. 
\end{re}

The following problem we found in the standard definition of gauge action was
the statement of the Algebraic Gauge-Invariant Uniqueness Theorem for Leavitt
path algebras. We can re-state it in the following form.

\begin{te}\label{AGIUT}
(The Schematic Algebraic Gauge-Invariant Uniqueness Theorem.)  Let $E$ be a graph,  $K$ any commutative unitary ring and  $A$ any $K$-algebra. 
Denote by $\rho\colon\GL_1\to\autn(L_K(E))$ the gauge-action of $L_K(E)$.
Assume that
$$\phi\colon L_K(E)\rightarrow A$$
 is a $K$-algebra homomorphism such that $\phi (r v)\ne 0$ for every $v\in E^0$ and $r\in K\setminus\{0\}$. If there exists an action
$\sigma \colon \GL_1 \rightarrow \autn(A)$ such that $(\phi\o 1)  \rho_R(z)=\sigma_R(z)(\phi\o 1)$ for every associative, commutative unital $K$-algebra $R$ and any $z\in R^\times$,  then $\phi $ is injective.
\end{te}

The proof is straightforward since by Remark \ref{notilla} the homomorphism $\phi$ is graded relative to the grading induced by $\sigma$ in $A$. Then we can apply \cite[Theorem 5.3, p. 476]{Tom_com}.\medskip

If $K$ is a field, we have the following:
\begin{co} Let $E$ be a graph,  $K$  any field and $A$ any $K$-algebra. 
Denote by $\rho\colon\GL_1\to\autn(L_K(E))$ the gauge-action of $L_K(E)$.
Assume that
$$\phi\colon L_K(E)\rightarrow A$$
 is a $K$-algebra homomorphism such that $\phi (v)\ne 0$ for every $v\in E^0$. 
 If there exists an action
$\sigma \colon \GL_1 \rightarrow \autn(A)$ such that $(\phi\o 1)  \rho_R(z)=\sigma_R(z)(\phi\o 1)$ for every associative, commutative unital $K$-algebra $R$ and any $z\in R^\times$,  then $\phi $ is injective.
\end{co}
\medskip

Also by using the gauge action in schematic sense the hypothesis on the infiniteness
of the ground field $K$ in \cite[Proposition 1.6]{Ahn} can be dropped. 
Since the notion of graded ideal and of gauge invariant ideal agree when we use the schematic version of the gauge action, such exceptionalities as the ones observed in \cite[Proposition 1.7]{Ahn} are no longer present. As pointed out previously the gauge action in schematic sense is more demanding. So it allows us to eliminate certain hypothesis simply because we are paying a higher price imposing conditions on all scalar extensions of the algebra.

\subsection{Cross product of algebras}
Let $K$ be a commutative unitary ring  and $A$ a $K$-algebra with an  action
$\rho\colon G \to\autn(A)$ where $G$ is an affine group scheme.
This means that $\rho$ is a natural transformation between the given
$K$-group functors. Then we define
\begin{de}
The fixed subalgebra $A^\rho$ of $A$ under $\rho$ is the one
whose elements are the elements $a\in A$ such that $\rho_R(z)(a\otimes 1)=
a\o 1$ for any $K$-algebra $R$ and any $z\in G(R)$.
\end{de}
If $A$ and $B$ are $K$-algebras provided with actions 
$\rho\colon G\to\autn(A)$ and $\sigma\colon G\to\autn(B)$
then there is an action $\rho\o\sigma\colon G\to\autn(A\o B)$ such that for
any $K$-algebra $R$ and any $z\in G(R)$ we have  
$(\rho\o\sigma)_R(z)$ given by the composition
\[
\xymatrix{
(A\o B)_R\ar@{.>}[d] _{(\rho\o\sigma)_{R}(z)} &=&  A\o B\o R\ar[r]^{1\o\delta} & A\o B\o R \o R\ar[r]^\theta &  A_R\o B_R\ar[d]^{\rho_R(z)\o\sigma_R(z^{-1})}\\
  (A\o B)_R & =& A\o B\o R & A\o B\o R\o R\ar[l]^{1\o\mu} & A_R\o B_R\ar[l]^{\theta^{-1}}}
\]
where: 
\begin{itemize}
\item $\delta\colon R\to R\o R$ is given by $\delta(z)=z\o 1$, 
\item $\theta$ is the isomorphism $a\o b\o r\o r'\mapsto a\o r\o b\o r'$,
\item $\mu\colon R\o R\to R$ is the multiplication $\mu(r\o r')=rr'$.
\end{itemize}
Summarizing $$(\rho\o\sigma)_R(z)=(1\o\mu)\theta^{-1}(\rho_R(z)\o\sigma_R(z^{-1}))\theta(1\o\delta).$$
Now a direct (but not short) computation reveals that $$(\rho\o\sigma)_R(zz')=
(\rho\o\sigma)_R(z)(\rho\o\sigma)_R(z')$$ for any $z$ and $z'$. So any $(\rho\o\sigma)_R(z)$ is invertible with inverse $(\rho\o\sigma)_R(z^{-1})$. Moreover, since
$(\rho\o\sigma)_R(z)$ is a composition of $R$-algebras homomorphisms, then
$(\rho\o\sigma)_R(z)\in\aut((A\o B)_R)$.
\begin{de}
The action $\rho\o\sigma\colon G\to\autn(A\otimes B)$ will be called the
tensor product action of $\rho$ and $\sigma$. The fixed point subalgebra
$(A\o B)^{\rho\o\sigma}$
of $A\o B$ under $\rho\o\sigma$ will be denoted $A\s{\rho}{\sigma}B$ and
called the cross  product of $A$ and $B$ by the actions $\rho$ and $\sigma$. 
If there is no ambiguity with respect to the actions involved we could shorten
the notation to $A\o_G B$.
\end{de}

Consider now two Leavitt path algebras $L_K(E)$ and $L_K(F)$ of the graphs $E$ and $F$ respectively. We assume given the gauge action of each algebra and ask about the cross product algebra  $L_K(E)\o_{\GL_1} L_K(F)$. 
With not much effort one can prove that it consists on the elements of the form
$\sum_{n\in\Z} a_n\o b_n$ where $a_n\in L_K(E)$ with $\deg(a)=n$ while
$b_n\in L_K(F)$ has also degree $n$.  Of course we can define on this algebra also an action $\tau\colon \GL_1\to\aut(L_K(E)\o_{\GL_1} L_K(F))$ by 
declaring for any $K$-algebra $R$ and each $z\in R^\times$ that 
$\tau_R(z)(a_n\otimes r\o b_n\o s)=a_n\o z^n r\o b_n \o s$. This action induces a grading on $L_K(E)\o_{\GL_1} L_K(F)$ in which the homogeneous component of degree $n$ is the $K$-submodule generated by the  elements of the form $a\o b$
where $a$ and $b$ are homogeneous of degree $n$.

Our next goal is to prove
\begin{te}
If $E$ and $F$ are row-finite graphs with no sinks, then
there is an isomorphism $L_K(E)\o_{\GL_1} L_K(F)\cong L_K(E\times F)$
where the product of the graph is the one described in Section 2.
\end{te}
Proof. 
For a graph $E$ denote by $\hat E$ the extended graph of $E$: the vertices
of $\hat E$ are those of $E$ and the arrows of $\hat E^1$ are those of $E^1$ plus  a family $\{f^*\colon f\in E^1\}$ of new edges such that $s(f^*)=r(f)$ and
$r(f^*)=s(f)$ for any $f\in E^1$. The path algebra $KE$ is the associative $K$-algebra with basis the set of all paths of $E$ (so it is free as a $K$-module).
There is a well known relation $L_K(E)\cong K\hat E/I$ where $I$ is the ideal
of $K\hat E$ generated by the CK relations.

We consider the path algebra $K(\widehat{E\times F})$, then there is a canonical homomorphism of $K$-algebras  map $K(\widehat{E\times F})\to  L_K(E)\o_{\GL_1} L_K(F)$ such that 
for any $(u,v)\in E^0\times F^0$ and $(f,g)\in E^1\times F^1$ we have
$$(u,v)\mapsto u\o v,\quad (f,g)\mapsto f\o g,\quad (f^*,g^*)\mapsto f^*\o g^*.$$
This homomorphism induces one $\phi\colon L_K(E\times F)\to L_K(E)\o_{\GL_1} L_K(F)$ such that $\phi(r(u,v))=r u\otimes v\ne 0$ for each $(u,v)\in E^0\times F^0$ and $r\in K\setminus\{0\}$ (see Remmark \ref{ffq}). Furthermore if we take the action $\tau\colon\GL_1\to\aut(L_K(E)\o_{\GL_1} L_K(F))$ defined above, we see that $(\phi\otimes 1)\rho_R(z)=
\tau_R(z)(\phi\o 1)$ where $\rho$ is the Gauge action of $L_K(E\times F)$.
Thus applying Theorem \ref{AGIUT} we conclude that $\phi$ is a monomorphism. To see that it is also an epimorphism we need the hypothesis that the graphs have no sinks. Since $L_K(E)\o_{\GL_1} L_K(F)$ is generated
by elements of the form $a\o b$ where $\deg(a)=\deg(b)$ it suffices to show that these elements are in the image of $\phi$.
First we prove that if $\mu$ and $\tau$ are paths of the same length (say $n$)
and $u$ is a  vertex, then $\mu\tau^*\o u$ is in the image of $\phi$:
Indeed, $\mu\tau^*\o u=\mu\tau^*\o\sum_i g_ig_i^*$ (since $F$ is row-finite and has no sink). If $\mu=f\mu'$ where $f\in E^1$ and $\mu'$ is a path then $\mu\tau^*\o u=\sum_i (f\o g_i)(\mu'\tau^*\o g_i^*)$ and if $\tau=h\tau'$ with $h\in E^1$ and $\tau'$ a path then
$\mu\tau^*\o u=\sum_i (f\o g_i)(\mu'\tau'^*h^*\o g_i^*)=\sum_i (f\o g_i)(\mu'\tau'^*\o r(g_i))(h^*\o g_i^*)$. Applying a
suitable induction hypothesis this proves that $\mu\tau^*\o u$ is in the image
of $\phi$. Symetrically it can be proved that the image of $\phi$ contains the
elements of the form $v\o \sigma\delta^*$ with $v\in E^0$ and $\sigma$, $\delta$ being paths of $F$ of the same degree.
Now any generator of $L_K(E)\o_{\GL_1} L_K(F)$ say $\mu\tau^*\o\sigma\delta^*$ such that $\deg(\mu)-\deg(\tau)=\deg(\sigma)-\deg(\delta)$ can be
written as a product of element which obey some of the followings patterns:
\begin{itemize}
\item $f\o g$ with $f\in E^1$ and $g\in F^1$.
\item $f^*\o g^*$ with $f\in E^1$ and $g\in F^1$.
\item $\mu\tau^*\o u$ with $u\in F^0$, $\mu$ and $\tau$ being paths of $E$ of the same length.
\item $v\o\sigma\delta^*$ with $v\in E^0$, $\sigma$ and $\delta$ being paths of $F$ of the same length.
\end{itemize}
Since any of these elements is in the image of $\phi$, this proves that $\phi$ is an epimorphism.\medskip

{\bf Acknowledgement.} The authors wishes to express their special thanks to the Referee for his/her careful reading and the many suggestions and corrections.

\end{document}